\documentclass[letterpaper,english,conference]{IEEEtran}
\usepackage[T1]{fontenc}
\usepackage{array}
\usepackage{float}
\usepackage{fancybox}
\usepackage{calc}
\usepackage{amsmath}
\usepackage{graphicx}

\makeatletter

\pdfpageheight\paperheight
\pdfpagewidth\paperwidth

\providecommand{\tabularnewline}{\\}
\floatstyle{ruled}
\newfloat{algorithm}{tbp}{loa}
\providecommand{\algorithmname}{Algorithm}
\floatname{algorithm}{\protect\algorithmname}

\usepackage{cite}

\usepackage{tikz}
\usetikzlibrary{decorations.pathmorphing}

\ifCLASSINFOpdf
   \graphicspath{{./images/}}
   \DeclareGraphicsExtensions{.pdf,.jpeg,.png}
\else
   \graphicspath{{../images/}}
   \DeclareGraphicsExtensions{.eps}
\fi
\usepackage{amssymb}
\usepackage{algorithm} 

\usepackage{multicol}
\usepackage{multirow}
\usepackage{booktabs}

\usepackage{hyperref}

\usepackage[tight,footnotesize]{subfigure}

\usepackage{stfloats}

\usepackage{url}

\usepackage{graphicx}

\usepackage{algpseudocode}

\makeatother

\usepackage{babel}
\begin{document}
\begin{center}
\textbf{\emph{}}%
\noindent\shadowbox{\begin{minipage}[t]{1\textwidth - 2\fboxsep - 2\fboxrule - \shadowsize}%
\begin{center}
\textbf{\emph{\LARGE{}This is a short (uncorrected) version of the
paper.}}
\par\end{center}{\LARGE \par}
\begin{center}
\textbf{\emph{\LARGE{}The final publication is available at Springer
via }}
\par\end{center}{\LARGE \par}
\begin{center}
\textbf{\emph{\LARGE{}https://link.springer.com/article/10.1007/s00521-017-2881-3 }}
\par\end{center}{\LARGE \par}
\begin{center}
\textbf{\emph{\LARGE{}DOI: 10.1007/s00521-017-2881-3}}
\par\end{center}{\LARGE \par}%
\end{minipage}}
\par\end{center}

\newpage{}

\title{Drone Squadron Optimization: a Self-adaptive Algorithm for Global
Numerical Optimization}

\author{\IEEEauthorblockN{Vin\'{i}cius V. de Melo\IEEEauthorrefmark{1} and
		Wolfgang Banzhaf\IEEEauthorrefmark{2}}
		
	\IEEEauthorblockA{\IEEEauthorrefmark{1}Institute of Science and Technology\\Federal University of S\~{a}o Paulo\\
		S\~{a}o Jos\'{e} dos Campos, SP, Brazil\\ Email: vinicius.melo@unifesp.br}

	\IEEEauthorblockA{\IEEEauthorrefmark{2}Department of Computer Science and Engineering and BEACON Center for the Study of Evolution in Action\\Michigan State University\\East Lansing, MI, 48864, USA\\ Email: banzhaf@msu.edu}
}

\maketitle
\begin{abstract}
This paper proposes Drone Squadron Optimization, a new self-adaptive
metaheuristic for global numerical optimization which is updated online
by a hyper-heuristic. DSO is an artifact-inspired technique, as opposed
to many algorithms used nowadays, which are nature-inspired. DSO is
very flexible because it is not related to behaviors or natural phenomena.
DSO has two core parts: the semi-autonomous drones that fly over a
landscape to explore, and the Command Center that processes the retrieved
data and updates the drones' firmware whenever necessary. The self-adaptive
aspect of DSO in this work is the perturbation/movement scheme, which
is the procedure used to generate target coordinates. This procedure
is evolved by the Command Center during the global optimization process
in order to adapt DSO to the search landscape. DSO was evaluated on
a set of widely employed benchmark functions. The statistical analysis
of the results shows that the proposed method is competitive with
the other methods in the comparison, the performance is promising,
but several future improvements are planned. 
\end{abstract}

\section{Introduction}

Meta-heuristics are general optimization methods used to solve a variety
of problems~\cite{osman2012meta}. When derivatives are costly or
impossible to obtain, derivative-free methods such as meta-heuristics
are usually employed for global optimization. Those methods make few
or no assumptions about the problem being optimized, being able to
deal with problems that are non-differentiable, irregular, noisy,
or dynamic. Otherwise, if derivatives are feasible, and one has a
convex problem, then convex optimization methods may be orders of
magnitude faster and more accurate.  Many meta-heuristic algorithms
are nature-inspired, with Genetic Algorithm (GA)~\cite{Goldberg1989}
being one of the most used. Another important category of meta-heuristic
is the Swarm Algorithm, introduced with the Ant Colony Optimization~\cite{dorigo1992ant}
for combinatorial optimization and the Particle Swarm Optimization
(PSO) algorithm in 1995~\cite{eberhart1995new} to solve continuous
optimization problems as in this work.

PSO has inspired the creation of many swarm intelligence algorithms,
such as Bees Algorithm, Artificial Bee Colony, Grey Wolf Optimizer,
Firefly algorithm, Bat algorithm, Gravitational Search Algorithm,
Glowworm swarm optimization, Cuckoo Search, among others. Fister Jr
et al.~\cite{fister2013brief} provide a short review and a list
of nature-inspired methods.

Here we change the source of inspiration. Instead of adopting a paradigm
from nature, we propose an \emph{artifact-inspired algorithm}\footnote{The terminology employed in this work is using the artifact as a metaphor
which - by way of analogy - can facilitate ones understanding.}, that is, it is inspired by something artificially created (human-made),
more specifically, \emph{drones}. In this paradigm, such an algorithm
is not bound by a particular realization at hand as drones are flexible
machines (at several aspects), not biological entities. Instead, it
can use a variety of different mechanisms/procedures without losing
its core characteristics. 

Another important aspect of the technique proposed here is that it
is self-adaptive regarding code modification, not only in the parameters
configuration. This means that the technique can change the procedures
that the drones use to generate solutions, meaning it can partially
evolve itself during the search. Some researchers did similar work,
but there are important differences that are shown later in this paper. 

Self-adaptation is one of the characteristics that can provide largest
improvements in performance~\cite{Back:1998:OPC:2379195.2379199,Qin_Suganthan_2005,AMALGAMvrugt2009self}.
However, most techniques use human-developed adaptation schemes that
do not cover every problem and may be unable to perform well on dynamic
situations. Thus, methods that can learn and self-adapt are of great
value. The key contributions of this paper are:
\begin{itemize}
\item the proposal of an artifact-inspired paradigm that is not tied to
any natural phenomena or behavior but could automatically act like
any of them;
\item a novel self-adaptive metaheuristic that can evolve itself on-the-fly,
transforming its behavior into that of other paradigms such as evolutionary
or swarm;
\item an explicit separation between the Controller and the semi-autonomous
exploration entities in a team approach.
\end{itemize}

The remaining of the paper is structured as follows: Section~\ref{sec:Drone-Squadron-Optimization}
presents our proposal. In Section~\ref{sec:Related-work} we present
related works. Section~\ref{sec:Experimental-results} demonstrates
the strength of the method with computational experiments. A further
discussion is given in Section~\ref{sec:Further-Discussion}. Section~\ref{sec:Conclusions-and-Future}
presents the conclusions and future work.

\section{\label{sec:Drone-Squadron-Optimization}Drone Squadron Optimization}

Drones, like the artificially built submarine machines or the well-known
flying machines, such as balloons, airplanes, helicopters, quad-copters,
can navigate autonomously or remotely, have sensors, can communicate
over large distances, can use solar power energy, and, one of the
most important features: can be upgraded/improved not only in terms
of hardware but also by changing their software~(the firmware). Therefore,
as the artificially built machine has a software~(firmware) to control
its behavior, researchers are free to add mechanisms to the algorithm
as common software upgrades, which is easier than looking for a natural
phenomena to justify the improvement.  

The Drone Squadron Optimization (DSO) proposed here may be related
to Particle Swarm Optimization, Artificial Bee Colony, or any other
Swarm algorithm, because it is based on the movement of entities in
search-space. However, as explained before, the movement of the squadron
is not necessarily based on behavior observed in nature. DSO's approach
allows it to automatically choose to use recombination and/or perturbation
of solutions with distinct procedures, making it act as an evolutionary
algorithm, swarm algorithm, probabilistic algorithm, or other, according
to how it performs on the search landscape. Moreover, those procedures
may have their \emph{actual code} updated during the search.

\begin{figure}[h]
\begin{centering}
\includegraphics[width=0.8\columnwidth]{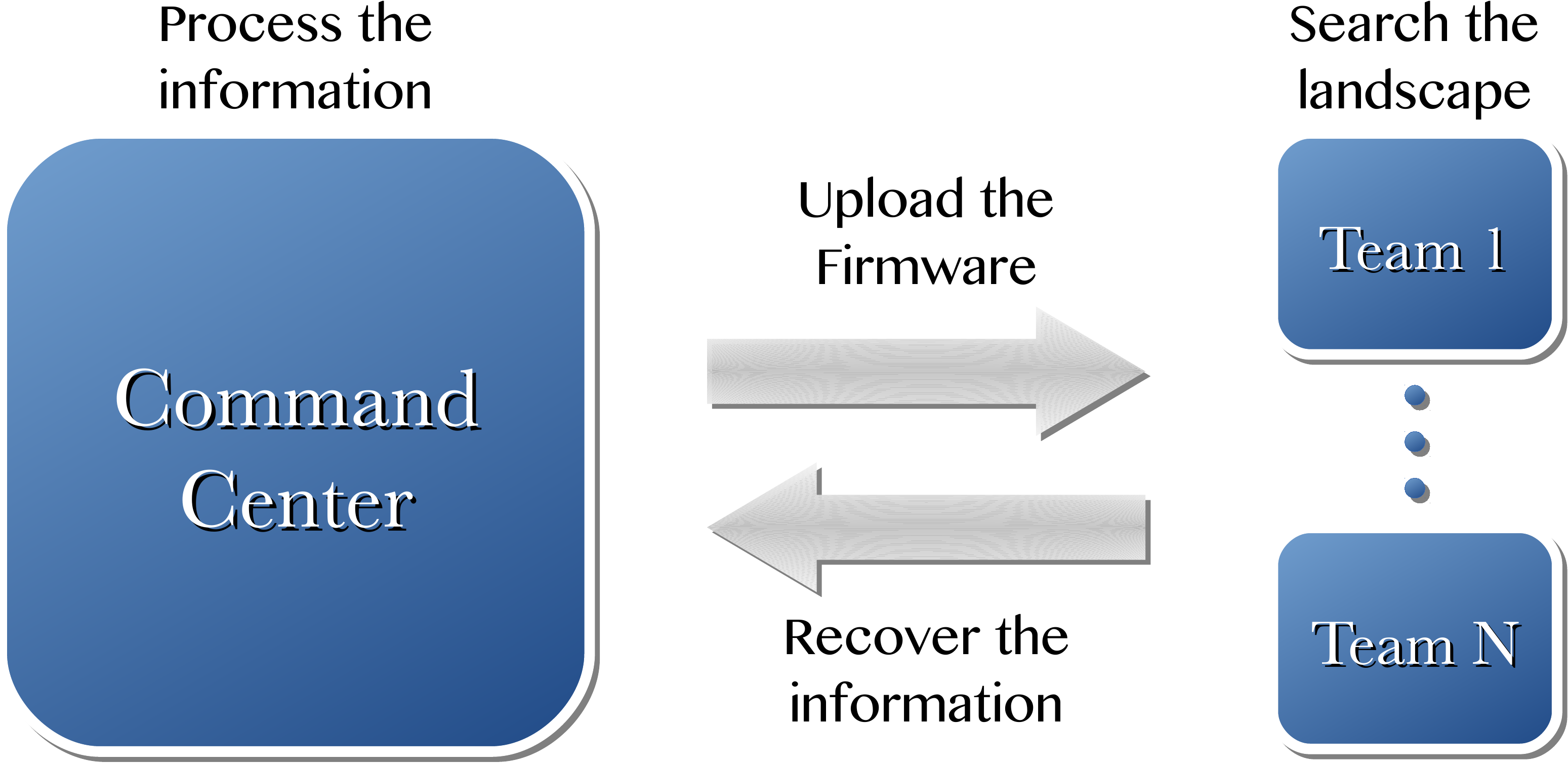}
\par\end{centering}
\caption{\label{fig:DSO-schema.}High-level DSO abstraction.}
\end{figure}

The DSO algorithm presented here is composed of a Drone Squadron with
different teams and a Command Center, which uses information collected
by/from the drones to maintain partial control of the search, and
to develop new firmware for controlling the drones (see Figure~\ref{fig:DSO-schema.}).
A drone \emph{is not} a solution; it moves to a coordinate which is
a solution. The firmware contains the procedures (codes) and configurations
used by the teams to search the landscape. In this work, the perturbation
procedure is an actual source code, a \emph{string} to be parsed and
executed by the drone.

Conceptually, the Command Center is a central place for carrying out
orders and for supervising tasks. It receives inputs, processes data,
and generates outputs from decisions made. The Command Center can
update the firmware of drones whenever it decides, dynamically adapting
drone behavior to the problem.

A group of drones is divided into teams of the same size (necessary
to the selection mechanism), where each team has its firmware that
controls the movement of each drone; there are, thus, distinct firmware
for each team. This means that \emph{each team has a distinct way
of sampling the search space from the same current set of solutions};
the teams are not like species nor niching in evolutionary algorithms.

The drones have a search mission on the landscape (the objective function)
to locate a particular target, whose value is obtained by the drone's
sensor. The teams are not necessarily multiple groups searching distinct
and distant regions of the landscape. In fact, all drones move from
specific departure points that may be the same for some of the teams,
but not for others. As the teams have distinct firmwares, even though
they move from the same departure points they may get to distinct
coordinates. Also, for the same reason teams may overlap search regions.
Moreover, a team does not have to follow another team, unless the
Command Center encode that in the firmware.

Algorithm~\ref{alg:high-level-dso} presents a high-level algorithm
of DSO, while the next subsections explain more details. One may observe
that DSO is more complex than traditional nature-inspired methods
as it has several other characteristics to deal with the Control Center,
the teams, the firmware adaptation, among others. 

\begin{algorithm*}[!tbh]
\caption{\label{alg:high-level-dso}High-level DSO algorithm. The abbreviations
are: $GlobalBestCoords$ (GBC), $\ensuremath{GlobalBestOFV}$ (GBOFV),
$\ensuremath{CurrentBestCoords}$ (CBC), $\ensuremath{CurrentBestOFV}$
(CBOFV), $TeamCoords$ (TmC), $TeamOFV$ (TmOFV), $\ensuremath{TrialCoords}$~(TC).}

\renewcommand{\algorithmicforall}{\textbf{for each}} 
\small
\begin{algorithmic}[1]
\State \textbf{Input:} objective function, problem bounds ($\overrightarrow{LB}$,$\overrightarrow{UB}$), user defined constants, number of teams ($t$), number of drones per team ($N$), maximum number of iterations ($MaxIt$), maximum number of iterations without improvement (MaxStagnation), options for the firmware
\State \textbf{Output:}  Best solution ($GBC$) and best Objective function value ($GBOFV$)\\

\State \textbf{Initialize} the Command Center, drones and teams in $TeamCoords$, where $TeamCoords$ are sets of teams containing the drones' coordinates
\State \textbf{Send} the drones to scan the landscape at ($t\times N$) random coordinates
\State \textbf{Select} the $N$ best coordinates to be $CBC$\\
\While {stopping criteria are not met}
	\State \Comment{Procedures performed by the Teams}
	\ForAll {$team \in [1,2,\dots,t]$}
		\ForAll {$drone \in [1,2,\dots,N]$}
			\State \textbf{Generate} $TC$ using the perturbation scheme of the current team
			\State \textbf{Choose} the recombination method and apply it to combine $TC$ with $CBC_{drone}$, resulting in the $TmC_{drone,team}$
			\State \textbf{Choose} out-of-boundary correction method and apply it to $TmC_{drone,team}$, saving the violations that occurred
			\State \textbf{Move} to new coordinate ($TmC_{drone,team}$), scan the landscape, and set $TmOFV_{drone,team}$ = objective function value
		\EndFor
	\EndFor\\
	\Comment{Procedures performed by the Command Center}
	\State \textbf{Select} the best coordinates found by the teams in the current iteration
	\State \textbf{Set} $improved$ = 0
	\ForAll {$drone \in [1,2,\dots,N]$}
		\Comment{Same drone id, but distinct teams}
		\State \textbf{Rank} $TmC_{drone,team}$ according to $TmOFV_{drone,team}$ 
		\State \textbf{Set} $bestIdx$ = index of the best $TmOFV_{drone,team}$
		\If { $TmOFV_{drone,bestIdx}$ < $CBOFV_{drone}$ }
			\State \textbf{Set} $CBOFV_{drone}$ = $TmOFV_{drone,bestIdx}$
			\State \textbf{Set} $CBC_{drone}$= $TmC_{drone,bestIdx}$
			\If { $TmOFV_{drone,bestIdx}$ < $GBOFV_{drone}$ }
				\State \textbf{Set} $improved = improved + 1$
				\State \textbf{Set} $GBOFV$ = $TmOFV_{drone,bestIdx}$
				\State \textbf{Set} $GBC$= $TmC_{drone,bestIdx}$
			\EndIf
		\EndIf
		\EndFor
	\If {$improved$ == 0} 
		\State \textbf{Check} and apply stagnation control
	\EndIf
	\State \textbf{Set} $TeamQuality$ as the average rank plus boundary violations
	\State \textbf{Update} the firmware, considering $TeamQuality$
\EndWhile
\end{algorithmic}
\end{algorithm*}

\subsection{\label{subsec:Perturbation}Drone movement}

The drones use an autonomous system to calculate target positions,
move to them, and collect the information that is sent back to the
Command Center. The following mechanisms available to DSO are employed
in various optimization techniques, either evolutionary or non-evolutionary.
Each mechanism may have more than one implementation giving DSO many
exploration and exploitation capabilities.

For the current DSO version, the firmware contains only the mechanism
to generate new trial coordinates (TC) through perturbation; thus,
it is the core of the drone's firmware. New firmware are generated
based on the well-known perturbation scheme of a biased random walk:
\begin{align}
P= & Departure+Offset(),\label{eq:Perturb}\\
TC & =calculate(P),
\end{align}

\noindent where \emph{$Departure$} is a coordinate (a solution point
in the search space), $Offset$ is a function that returns the perturbation
movement (a numerical value), and $P$ is the complete perturbation
formula that has to be calculated to return the trial coordinates\footnote{\noindent It is important to notice that all solutions are one-dimension
arrays; therefore, all operations present in this work are element-wise. }. This pattern of moving from a \emph{$Departure$} point is important
because it avoids the algorithm to simply shrink the search space
towards the origin\footnote{Many well-known benchmark functions have their global optimum at the
origin and there are algorithms that exploit such characteristic to
achieve high performance.}. In fact, even though the $Offset$ is shrinking to zero, it is still
a noise applied to the $Departure$ coordinates, which results in
a neighborhood search.

Instead of using a hard-coded procedure, $P$ is automatically updated
during the global optimization process and adapted to the problem.
The adaptation will be explained further as it is not performed by
the drone itself. Both the $Departure$ coordinates and the function
that gives the $Offset$ are modified, that is, distinct teams may
have distinct ways of choosing departure coordinates and how to calculate
the offset. As examples, suppose the two following perturbations for
teams~1 and~2, where $C_{1}$ is a user-defined constant, $G(0,1)$
is a scalar sampled from a Gaussian distribution with zero mean and
unit standard deviation, and $U(0,1,D)$ is an array of \emph{$D$}
columns of numbers sampled from a Uniform distribution with minimum~0
and maximum~1. The meaning of the terms can be seen in Algorithm~\ref{alg:high-level-dso}.

\begin{flushleft}
\begin{equation}
P_{1}:GBC+(C_{1}\times(GBC-CBC_{drone})).\label{eq:Perturb-Ex1}
\end{equation}
\par\end{flushleft}

\begin{flushleft}
\begin{equation}
P_{2}:CBC_{drone}+(G(0,1)\times(\sqrt{U(0,1,D)}+CBC_{drone})).\label{eq:Perturb-Ex2-1}
\end{equation}
\par\end{flushleft}

As one may observe, the two $Departures$ are the arrays $GBC$ and
$CBC_{drone}$, while the added expressions are the $Offsets$. The
perturbation scheme in DSO was implemented using a \emph{tree-structure}
representation with terminal and non-terminal nodes, and always follows
the pattern shown in Equation~\ref{eq:Perturb}, that is, a sum of
two terms where the first one is selected from a particular subset
of the terminals, and the second one is an expression grown using
the available terminals and non-terminals. 

Initially, the perturbation schemes were randomly generated (arbitrary
equations). However, after preliminary experiments it was identified
that a completely random set of initial perturbation schemes shows
poor performance. To overcome this issue, a set of \emph{reference
perturbations} was defined, that is, equations that will be used as
\emph{initial} perturbations for the teams, but that may be replaced
during the optimization process.

When exploiting a particular region of the landscape, the drones may
generate identical or almost identical target coordinates. Convergence
avoidance mechanisms are used to help the drones escape from local
optima. The use of scaling (Gaussian, Uniform, etc.), stagnation detection
and its treatment, and generation of opposite coordinates allows movements
to regions far from the neighborhood. Unfortunately, this may also
slow down reaching the optimum solution and/or reduce its accuracy
when a small budget is considered, but we chose to expand the exploration
capability of the method to locate promising regions in the search
space.

After the perturbation step generates the trial coordinates for each
drone ($TC_{drone}$), either a recombination to generate $TmC$ is
performed with the best coordinates found so far or no recombination
is done ($TmC=TC$). Currently, this selection is random and all recombination
procedures available to the drones have the same probability of selection.
Also, recombination is performed after perturbation, but changing
the order is a perfectly plausible option. That will change the behavior
of the method without invalidating the original artifact inspiration.

Finally, the drones may be allowed to move only inside a particular
perimeter. Therefore, if coordinates in $TmC$ are outside such perimeter,
then a correction must be made. As in the recombination, some correction
procedures are available to be chosen, and there is no bias in the
current DSO.

\subsection{Command Center}

By using the information collected by the drones, the Command Center
develops and updates new firmware to upload to them. The drones use
the new firmware and all data distributed by the Command Center to
control their behavior. Conceptually, cheap data analysis may be performed
on the drones due to their limited computing power, whereas the expensive
ones are done on the Command Center and shared with the drones.

\subsubsection{Firmware update}

The Command Center uses two information to measure the performance
of a team: 1) its rank regarding the objective function value, and
2) the degree of out-of-bound coordinates that were generated. It
is important to consider the violations because good solutions may
be generated by the correction procedures just by chance, while \emph{$TmC$}
had, in fact, big violations.

Given that each team generates $N$ solutions, and that the ID of
a drone in a \emph{team} is it index $drone=1,...,N$, the rank for
each ID is obtained by sorting $TmOFV_{drone,team}$ of distinct teams.
For an example with three teams, sort $TmOFV_{drone,1}$, $TmOFV_{drone,2}$,
$TmOFV_{drone,3}$. 

Regarding the violations, they are calculated as:

{\scriptsize{}}{\scriptsize \par}

{\small{}
\begin{equation}
violation_{drone,team}=\sum_{j=1}^{D}\begin{cases}
|TmC_{drone,team,j}-UB_{j}|\\
+\\
|LB_{j}-TmC_{drone,team,j}|
\end{cases},\label{eq:Violations}
\end{equation}
}{\small \par}

\noindent where \emph{D} is the problem's dimension, and considering
only the cases where $TmC_{drone,team,j}>UB_{j}$ or $TmC_{drone,team,j}<LB_{j}$.
The rank is calculated at every iteration and can be accumulated or
averaged after a series of iterations. 

As soon as a firmware updating criterion is reached (amount of iterations,
for instance), the Command Center replaces the \emph{w} worst firmware
by variations of the \emph{w} best firmware, i.e., for $w=1$ the
team with worst accumulated rank has its firmware updated with a variant
(random sub-tree replacement) of the best team firmware. There is
no recombination of codes. The new variant has to satisfy the following
rules: 
\begin{enumerate}
\item The size \emph{S} of the new perturbation \emph{$P_{k}$}, where $S$
is the number of nodes in the tree data structure, and $k$ is the
index of the worst of the $t$ teams, has to be $S(P_{k})>s_{min}$
and $S(P_{k})<s_{max}$, where $s_{min}$ and $s_{max}$ are user
defined parameters; 
\item The new perturbation has to be distinct from the original one, but
the current version only detects syntactic differences, not semantic
ones;
\item A function is not allowed to receive the same argument for the two
parameters, for instance, $sub(Shift,\,Shift)$;
\item The \emph{w} reference perturbations must not be replaced, thus they
are \emph{fixed perturbations};
\item The perturbation scheme of Equation~\ref{eq:Perturb} must hold.
\end{enumerate}
Consequently, one expects performance improvement after replacing
the firmware that gave the worst results by a variation of the firmware
that gave the best results. With respect to the \emph{fixed} perturbations,
this mechanism is to provide at least one firmware with reasonable
search capability. A \emph{fixed} perturbation may be one that favors
exploration, while the remaining one may be free to perform exploitation.

\subsubsection{Selection for next iteration, stagnation detection and treatment }

After all drones return their findings on the landscape (the objective
function values), the Command Center must decide which information
is important to keep in the search plan in order to be used in the
future. A hard selection mechanism (as shown in Algorithm~\ref{alg:high-level-dso})
chooses the best between the current solution and the new one, considering
the objective function values. However, stagnation must be treated. 

Stagnation is detected when the objective function value of the current
best solution at a particular index remains the same after a certain
number of iterations. When detected, a mechanism to allow further
exploration is used and soft selection is performed:

~

\textbf{\footnotesize{}if}{\footnotesize{} $(\ensuremath{TmOFV_{drone,bestIdx}<\ensuremath{CBOFV_{drone}}}$
}\textbf{\footnotesize{}or }{\footnotesize{}$U(0,1)<P_{acc})$}{\footnotesize \par}

\textbf{\footnotesize{}then}{\footnotesize{} $CBC_{drone}\ensuremath{=}TmC_{drone,bestIdx}$, }{\footnotesize \par}

~

\noindent where $bestIdx$ is the index of the best drone (rank 1
considering all teams), and $P_{acc}$ is the probability of accepting
a worse solution. Thus, some coordinates that resulted in lower-quality
solutions are inserted in the search plan, replacing higher-quality
solutions. Elitism concept is applied here to keep the best solution
in $CBC$.

\section{\label{sec:Related-work}Related work}

Here we consider only the most closely related work, i.e. work that
employs some mechanism to evolve a meta-heuristic regarding the mutation/perturbation
operator. One may divide the research into two basic evolutionary
branches: offline and online. In offline evolution, the main method
is a hyper-heuristic~\cite{Burke:2010:HBMH}, such as a GP-like method~\cite{banzhaf1998genetic},
responsible for generating new code for the meta-heuristic and running
it for some repetitions on one or more benchmark problems. The statistical
result of the runs is returned to the hyper-heuristic as the fitness
value. In online evolution, on the other hand, the hyper-heuristic
and the meta-heuristic work in co-operation, which means that there
is a single meta-heuristic run, not a set of repetitions, and the
meta-heuristic's code is evolved during the optimization by the hyper-heuristic.
Therefore, while in the offline approach the hyper-heuristic tries
to find the best static code to solve a single problem or various
problems, in the online approach the code is evolved for the particular
problem in a particular run. Given that DSO is an online approach,
below we consider only similar research from the literature. Unfortunately,
we could find only a few work.

Rashid and Baig~\cite{PSOGP} proposed PSOGP, an extension of~\cite{Poli:2005:EPS:2136684.2136712}
that works online, evolving the code during the optimization. Each
particle in the population has its force generating function, which
is evolved by GP. This approach gave PSO better exploration abilities
by slowing convergence and increasing chances of escaping from local
optima. 

In a similar approach, Si et al.~\cite{si2014grammatical} replaced
GP by GDE (Grammatical Differential Evolution~\cite{GDE2006grammatical})
and also achieved better performance than the PSO version used in
comparison.

One may also consider as online approaches some adaptive methods such
as those proposed in~\cite{PAPpeng2010population,AMALGAMvrugt2009self,SADEqin2009differential,whitacre2006credit},
where several methods or operators are coded and available for selection
during the optimization process. A mechanism, usually a probability
of success in solution improvement, is used to select the most promising
codes. However, the options are hard-coded, that is, the code is immutable.

\section{\label{sec:Experimental-results}Experimental results}

This next section presents a comparison of Drone Squadron Optimization
with state-of-the-art evolutionary and swarm optimization methods
in solving a set of well-known unconstrained/box-constrained continuous
global optimization benchmark functions.

\subsection{Benchmark functions}

DSO was tested on functions $f1$ to $f10$ of the benchmark problem
set from CEC' 2005 Special Session in Real Parameter Optimization~\cite{suganthan2005problem}.
Those problems present a diverse set of features such as multimodality,
ruggedness, ill-conditioning, interdependency, among others. In our
experiments, all functions were solved in 10 dimensions, with a maximum
\emph{total} of 100,000 objective function evaluations \emph{shared}
by the teams (it is not per team), and all configurations suggested
in~\cite{suganthan2005problem}.

\subsection{Implementation and configuration of the algorithms}


The experiments were executed on an Intel(R) Xeon(R) X5550 @ 2.67
GHz with 8GB RAM, Linux enterprise 3.14.1-gentoo i686. DSO was implemented
in MATLAB R2012b. The information and functions available to the Command
Center to generate perturbation schemes are shown below.
\begin{itemize}
\item \textbf{Departure coordinates:} $CurrentBestCoords$, $PermutedCurrentBestCoords$,
$PBestCurrentBestCoords$, Multivariate normal sampling (MVNS): new
random solutions sampled using the average and covariance matrix of
the $p$-$Best$ solutions found, $Opposition(CurrentBestCoords)$
as in~\cite{journals/tec/RahnamayanTS08}.
\item \textbf{Other terminals used in $Offset$:} $matInterval_{D}=(\overrightarrow{UB_{D}}-\overrightarrow{LB_{D}})$,
C1, C2, C3, $U(0,1)$, $U(0.5,1)$, $G(0,1)$, $abs(G(0.5,0.1))$,
$abs(G(0,0.01))$, $Step(CBC)=\sigma*G(0,1)_{N,D}*matInterval_{D}*U(0,\,0.5)$
as used in~\cite{de2014modified}, $Shift=TmC-CBC$, and $GBC$.
\item \textbf{Functions:} \textit{abs, neg, plus}, \textit{times}, \textit{sub},
\textit{protected division}, \textit{average}\emph{\_two\_args.}
\end{itemize}
DSO configuration is shown in Table~\ref{tab:ConfigDSO-Exp2} and
was empirically chosen after a few runs with distinct configurations,
but we intend to do a deeper investigation on the influence of the
parameters in a future work. $rand/1$ is the classical mutation strategy
from the Differential Evolution (but not linked to a particular crossover),
and MVNS+Step is \emph{inspired} by the CMA-ES technique, but employing
only sample generation and step calculation with $\sigma=0.04\times\mu_{eff}\times||\mu||$.
This formula is from the CMA-ES author's source code and was not tuned
to be used in DSO.

DSO will work with 4 teams of 25 members, giving the overall sample
size 100. The sample is not necessarily 100 distinct solutions that
will be perturbed to give new 100 solutions as performed by other
methods because the $Departure$ coordinates may be the same for all
teams at some iteration.

\begin{table}[tbh]
\caption{\label{tab:ConfigDSO-Exp2}Run and evolutionary parameter values for
DSO to solve CEC'2005 problems.}

\noindent \centering{}%
\begin{tabular}[t]{>{\centering}p{0.35\columnwidth}|>{\centering}p{0.45\columnwidth}}
\hline 
\textbf{Parameter}  & \textbf{Value}\tabularnewline
\hline 
\hline 
C1, C2, C3  & 0.5, 0.4, 0.9\tabularnewline
Teams  & 4\tabularnewline
Reference perturbations & \emph{rand}/1 using C1 as F, MVNS+Step\tabularnewline
Firmware update  & every iteration\tabularnewline
$w$ ($\#$ firmware updated)  & 1\tabularnewline
Fixed perturbation & \emph{rand}/1 using C1 as F\tabularnewline
Recombination & No recombination, Binomial recombination~\cite{book:price2005},
Exponential recombination~\cite{book:price2005}\tabularnewline
Stagnation  & 50 iterations\tabularnewline
Elitism when stagnated  & 1\tabularnewline
$P_{acc}$ & $10\%$ \tabularnewline
tree-size ($s_{min}$, $s_{max}$)  & 5, 20\tabularnewline
$p$-$Best$  & $25\%$ \tabularnewline
$CR$  & $U(0.4,0.9)$\tabularnewline
\hline 
\end{tabular}
\end{table}

\subsection{Performance comparison}

Results are evaluated over 36 runs (instead of 25 as the other methods)
because the method was executed on a cluster and all available computing
nodes were used. A successful run presents error lower than 1E-09.

The rank analysis is first done by Friedman rank sum test. Then, if
the statistical difference among the methods is detected, one employed
the Pairwise posthoc Test for Multiple Comparisons of Mean Rank Sums
for Unreplicated Blocked Data (Nemenyi-Test) to identify differences
between DSO (as the control method) and each method used in the comparison.
We assumed a significance level $\alpha=0.05$. 

\subsection{Results and Discussion}

In this section, we evaluate DSO on the first ten benchmark functions
of the problem set from CEC' 2005 Special Session in Real Parameter
Optimization\footnote{We intend to test all 25 benchmark functions and other problem sizes
in a future work.}. We focused on the dimension $D=10$ as investigated in~\cite{garcia2009study}.

Descriptive statistics of DSO runs are shown in Table~\ref{tab:Statistics_DSO}.
DSO failed for functions f8 and f10 while showed poor performance
for functions f3 and f7. The minimum median achieved values suggest
that DSO was unable to perform local optimization for these functions
as it could not reach the desired precision. Also, mean values are
much higher than median values for several functions, meaning outlier
failures. This behavior was not expected given the stagnation and
escaping procedures employed in the method. Hence, more investigation
is necessary to improve this aspect of the method.

\begin{table}[!h]
\caption{\label{tab:Statistics_DSO}Descriptive statistics of DSO runs obtained
in CEC\textquoteright 2005 Special Session in dimension $D=10$.}

\noindent \centering{}\resizebox{1.0\columnwidth}{!}{%
\begin{tabular}{c|cccccc}
\hline 
\textbf{Fun} & \textbf{Min} & \textbf{Median} & \textbf{Max} & \textbf{Mean} & \textbf{Std.Dev} & \textbf{SR}\tabularnewline
\hline 
\hline 
\textbf{f1} & 0.00E+00 & 0.00E+00 & 4.08E-19 & 1.13E-20 & 6.79E-20 & 1\tabularnewline
\textbf{f2} & 7.09E-21 & 4.06E-16 & 6.02E-14 & 4.03E-15 & 1.19E-14 & 1\tabularnewline
\textbf{f3} & 1.90E-07 & 4.99E-02 & 2.41E+03 & 7.70E+01 & 4.03E+02 & 0.028\tabularnewline
\textbf{f4} & 9.39E-17 & 5.24E-14 & 1.51E-10 & 5.72E-12 & 2.57E-11 & 1\tabularnewline
\textbf{f5} & 0.00E+00 & 0.00E+00 & 2.15E-01 & 5.97E-03 & 3.58E-02 & 0.944\tabularnewline
\textbf{f6} & 1.04E-19 & 4.47E-12 & 3.99E+00 & 3.32E-01 & 1.12E+00 & 0.917\tabularnewline
\textbf{f7} & 7.40E-03 & 5.66E-02 & 1.57E-01 & 6.53E-02 & 3.43E-02 & 0.028\tabularnewline
\textbf{f8} & 2.00E+01 & 2.00E+01 & 2.05E+01 & 2.01E+01 & 1.06E-01 & 0\tabularnewline
\textbf{f9} & 0.00E+00 & 0.00E+00 & 1.99E+00 & 4.70E-01 & 6.93E-01 & 0.639\tabularnewline
\textbf{f10} & 2.98E+00 & 1.34E+01 & 2.79E+01 & 1.37E+01 & 6.16E+00 & 0\tabularnewline
\hline 
\end{tabular}}
\end{table}

The results shown in Table~\ref{tab:Comparison-CEC05} are compared
with those from Table 13 in~\cite{garcia2009study}, which considers
only the average error for each problem. The optimization methods
in the comparison are: Hybrid Real-Coded Genetic Algorithm with Female
and Male Differentiation (BLX-GL50~\cite{garcia2005hybrid}), Real-Coded
Memetic Algorithm (BLX-MA~\cite{molina2005adaptive}), Cooperative
Evolution EA (CoEVO~\cite{povsik2005real}), canonical Differential
Evolution (DE~\cite{ronkkonen2005real}), Dynamic multi-swarm particle
swarm optimizer with local search (DMS-L-PSO~\cite{liang2005dynamic}),
Estimation of Distribution Algorithm (EDA~\cite{yuan2005experimental}),
Covariance Matrix Evolution Strategy and Restarting method (G-CMA-ES~\cite{auger2005restart}),
Steady-State Evolutionary Algorithm (K-PCX~\cite{sinha2005population}),
Covariance Matrix Evolution Strategy Improved with Local Search (L-CMA-ES~\cite{auger2005performance}),
Self-adaptive differential evolution algorithm for numerical optimization
(L-SaDE~\cite{Qin_Suganthan_2005}), and Steady-State Genetic Algorithm
(SPC-PNX~\cite{ballester2005real}).

\begin{table}[!t]
\caption{\label{tab:Comparison-CEC05}Average error obtained in CEC\textquoteright 2005
Special Session in dimension $D=10$. Results \textbf{\emph{worse}}
than those of DSO are in italic and bold-face.}

\noindent \centering{}\resizebox{1.0\columnwidth}{!}{%
\begin{tabular}{c|ccccc}
\hline 
\textbf{Algorithm} & \textbf{f1} & \textbf{f2} & \textbf{f3} & \textbf{f4} & \textbf{f5}\tabularnewline
\hline 
\hline 
\textbf{DSO} & 1.000E-09 & 1.000E-09 & 7.697E+01 & 1.000E-09 & 5.968E-03\tabularnewline
\textbf{BLX-GL50} & 1.000E-09 & 1.000E-09 & \textbf{\emph{5.705E+02}} & 1.000E-09 & 1.000E-09\tabularnewline
\textbf{BLX-MA} & 1.000E-09 & 1.000E-09 & \textbf{\emph{4.771E+04}} & \textbf{\emph{1.997E-08}} & \textbf{\emph{2.124E-02}}\tabularnewline
\textbf{CoEVO} & 1.000E-09 & 1.000E-09 & 1.000E-09 & 1.000E-09 & \textbf{\emph{2.133E+00}}\tabularnewline
\textbf{DE} & 1.000E-09 & 1.000E-09 & 1.940E-06 & 1.000E-09 & 1.000E-09\tabularnewline
\textbf{DMS-L-PSO} & 1.000E-09 & 1.000E-09 & 1.000E-09 & \textbf{\emph{1.885E-03}} & 1.138E-06\tabularnewline
\textbf{EDA} & 1.000E-09 & 1.000E-09 & 2.121E+01 & 1.000E-09 & 1.000E-09\tabularnewline
\textbf{G-CMA-ES} & 1.000E-09 & 1.000E-09 & 1.000E-09 & 1.000E-09 & 1.000E-09\tabularnewline
\textbf{K-PCX} & 1.000E-09 & 1.000E-09 & 4.150E-01 & \textbf{\emph{7.940E-07}} & \textbf{\emph{4.850E+01}}\tabularnewline
\textbf{L-CMA-ES} & 1.000E-09 & 1.000E-09 & 1.000E-09 & \textbf{\emph{1.760E+06}} & 1.000E-09\tabularnewline
\textbf{L-SaDE} & 1.000E-09 & 1.000E-09 & 1.672E-02 & \textbf{\emph{1.418E-05}} & \textbf{\emph{1.200E-02}}\tabularnewline
\textbf{SPC-PNX} & 1.000E-09 & 1.000E-09 & \textbf{\emph{1.081E+05}} & 1.000E-09 & 1.000E-09\tabularnewline
\hline 
\textbf{Algorithm} & \textbf{f6} & \textbf{f7} & \textbf{f8} & \textbf{f9} & \textbf{f10}\tabularnewline
\hline 
\textbf{DSO} & 3.322E-01 & 6.526E-02 & 2.005E+01 & 4.698E-01 & 1.365E+01\tabularnewline
\textbf{BLX-GL50} & 1.000E-09 & 1.172E-02 & \textbf{\emph{2.035E+01}} & \textbf{\emph{1.154E+00}} & 4.975E+00\tabularnewline
\textbf{BLX-MA} & \textbf{\emph{1.490E+00}} & \textbf{\emph{1.971E-01}} & \textbf{\emph{2.019E+01}} & 4.379E-01 & 5.643E+00\tabularnewline
\textbf{CoEVO} & \textbf{\emph{1.246E+01}} & 3.705E-02 & \textbf{\emph{2.027E+01}} & \textbf{\emph{1.919E+01}} & \textbf{\emph{2.677E+01}}\tabularnewline
\textbf{DE} & 1.590E-01 & \textbf{\emph{1.460E-01}} & \textbf{\emph{2.040E+01}} & \textbf{\emph{9.550E-01}} & 1.250E+01\tabularnewline
\textbf{DMS-L-PSO} & 6.892E-08 & 4.519E-02 & 2.000E+01 & 1.000E-09 & 3.622E+00\tabularnewline
\textbf{EDA} & 4.182E-02 & \textbf{\emph{4.205E-01}} & \textbf{\emph{2.034E+01}} & \textbf{\emph{5.418E+00}} & 5.289E+00\tabularnewline
\textbf{G-CMA-ES} & 1.000E-09 & 1.000E-09 & 2.000E+01 & 2.390E-01 & 7.960E-02\tabularnewline
\textbf{K-PCX} & 4.780E-01 & \textbf{\emph{2.310E-01}} & 2.000E+01 & 1.190E-01 & 2.390E-01\tabularnewline
\textbf{L-CMA-ES} & 1.000E-09 & 1.000E-09 & 2.000E+01 & \textbf{\emph{4.490E+01}} & \textbf{\emph{4.080E+01}}\tabularnewline
\textbf{L-SaDE} & 1.199E-08 & 2.000E-02 & 2.000E+01 & 1.000E-09 & 4.969E+00\tabularnewline
\textbf{SPC-PNX} & \textbf{\emph{1.891E+01}} & \textbf{\emph{8.261E-02}} & \textbf{\emph{2.099E+01}} & \textbf{\emph{4.020E+00}} & 7.304E+00\tabularnewline
\hline 
\end{tabular}}
\end{table}

Friedman rank sum test was performed to compare the 12 techniques'
averages. As the Friedman test indicates significant difference ($\chi^{2}$
(11) = 20.6765, $p=3.688E\textrm{{-}}02<\alpha$), it is meaningful
to conduct multiple comparisons in order to identify differences between
the techniques on the 10 problems. For such purpose, Table~\ref{tab:Comparisons-WTL-Pvalues}
presents the comparisons regarding \emph{win}, \emph{tie}, and \emph{loss},
as well as Nemenyi's posthoc \emph{p}-values.

On the one hand, with respect to the W/T/L values, DSO was better
than BLX-MA, CoEVO, EDA, K-PCX, and SPC-PNX. On the other hand, as
can be noticed by the \emph{p}-values, no significant differences
($p<\alpha$) were detected by the posthoc analysis for $\alpha=0.05$,
even though G-CMA-ES and DMS-L-PSO won 7 out of 10 times (lower average
error), and DSO won 6 times when compared to BLX-MA. Therefore, one
may conclude that DSO is as good as the methods selected for that
Special Session, although it is based on an construction mechanism
which sometimes generates perturbation schemes that give low quality
solutions and waste precious evaluations, while other methods are
hand-made specialized developments. Thus, this characteristic is a
good opportunity of improving DSO by trying to learn to avoid such
poor schemes.

Processing time cannot be compared as not all methods have source
code freely available online. Nevertheless, DSO was developed in Matlab
and requires parsing and executing the perturbation source code, which
gives a much longer running time. 

\begin{table*}[htpb]
\caption{\label{tab:Comparisons-WTL-Pvalues}Comparison of the averages: DSO
versus other methods (Wins/Ties/Losses). \emph{p}-value was calculated
by Nemenyi's post-hoc test.}

\noindent \centering{}%
\begin{tabular}{ccccccc}
\hline 
 & \textbf{BLX-GL50} & \textbf{BLX-MA} & \textbf{CoEVO} & \textbf{DE} & \textbf{DMS-L-PSO} & \textbf{EDA}\tabularnewline
\hline 
\hline 
\textbf{W/T/L} & 3/3/4 & 6/2/2 & 5/3/2 & 3/3/4 & 1/2/7 & 3/3/4\tabularnewline
\textbf{\emph{p}}\textbf{-value} & 1.00 & 1.00 & 1.00 & 1.00 & 0.99 & 1.00\tabularnewline
\hline 
 & \textbf{G-CMA-ES} & \textbf{K-PCX} & \textbf{L-CMA-ES} & \textbf{L-SaDE} & \textbf{SPC-PNX} & \tabularnewline
\hline 
\hline 
\textbf{W/T/L} & 0/3/7 & 4/2/4 & 3/2/5 & 2/2/6 & 5/3/2 & \tabularnewline
\textbf{\emph{p}}\textbf{-value} & 0.50 & 1.00 & 1.00 & 1.00 & 1.00 & \tabularnewline
\hline 
\end{tabular}
\end{table*}

\section{\label{sec:Further-Discussion}Further Discussion}

In this section, we discuss some aspects of DSO that could raise questions
and then clarify them. First of all, DSO employs ideas of some distinct
frameworks, thus it does not completely fits to a particular one.

DSO can adapt its parameters to produce better performance during
the optimization process, as similar adaptive techniques~\cite{DBLP:journals/tec/ZhangS09,Qin_Suganthan_2005,thierens2005adaptive,whitacre2006credit}.
However, DSO does more than that.

In DSO, the perturbation operator is also modified, not just the parameters.
In this case, DSO can be seen as a generative hyper-heuristic approach~\cite{Burke:2010:HBMH},
where code is evolved instead of parameters. However, in a standard
hyper-heuristic application, the generation of codes is the main procedure,
while in DSO it is a component of the method. Thus, DSO is a hybrid
approach with co-evolution of solutions to the problem being optimized
and the procedures to optimize it.

DSO has independent teams operating on current solutions, which is
similar to distributed EAs (dEAs)~\cite{lin1994coarse,alba2002parallelism,tasoulis2004parallel}.
However, most dEAs employ the same algorithm on all subpopulations
(either with the same or different control parameters). On the other
hand, in DSO distinct teams may use the same $Departure$ coordinates
(same set of solutions - population), but with distinct ways of calculating
the $Offset$. At other period of the search, the teams may have the
same $Offset$ formulation, but distinct $Departures$. Given that
the teams may have distinct perturbation and recombination procedures,
resulting in distinguishable behaviors, they could be seen as distinct
algorithms.

In dEAs, the parallel algorithms have their own populations searching
distinct regions of the search-space, and a migration procedure is
adopted to try to escape from local optima and to speed-up search.
On the other hand, the teams in DSO are not sub-populations; hence
no migration is required. Another important argument why migration
is not used is because DSO's teams are neither directly co-evolving
nor cooperating. In fact, the selection mechanism considers the solutions
from all teams; therefore they are primarily competing to each other,
even though all teams may improve the same shared information. Nevertheless,
a comparison between competition and cooperation is out of the scope
of this paper.

By using distinct algorithms at the same time, DSO could be a portfolio~\cite{PAPpeng2010population}
of evolving (adapting) evolutionary processes. Most differences for
dEAs can be applied to the portfolio framework, but in portfolios
and other similar methods, such as AMALGAM~\cite{AMALGAMvrugt2009self}
and A-Teams~\cite{rachlin1999teams}, distinct populations use distinct
algorithms. However, again, these algorithms are hard-coded and are
not improved during the optimization.

Finally, it is important to recall that DSO's teams are semi-autonomous;
they must obey the Command Center. The Command Center can upgrade
the firmwares, select the solutions found by the teams, change the
amount of teams and drones in each team, restart the solutions, perform
analyzes of the solutions, start a local-search, among other possibilities.

Regarding computation time, the current DSO implemented in Matlab
is certainly slower than all methods used in comparison due to the
fact it generates code during the optimization process, and that new
code has to be evaluated (parsed and compiled) to give the desired
results. Once compiled, the perturbation scheme is as fast as hard-coded
functions. However, if the perturbation scheme is often changed, then
a slowdown is clearly noticeable.

\section{\label{sec:Conclusions-and-Future}Conclusions and Future work}

The Drone Squadron Optimization (DSO) introduced in this paper is
a novel self-adaptive method for global numerical optimization. It
is a population-based search algorithm, but not nature-inspired. Because
it is artifact-inspired, it was developed to be a very flexible technique,
allowing the insertion and removal of components without disrupting
the original concept. For instance, DSO may or may not apply recombination,
and any kind of data recombination is permitted as there is no biology
concept involved. DSO is self-adaptive from its first version, because
this is the role played by the Command Center.

DSO was evaluated on many well-known continuous benchmark problems
and compared to many state-of-the-art methods. After the hypothesis
test, no significant difference was found, suggesting that DSO showed
competitive performance. However, it is important to remember that
DSO is a technique that evolves itself during the optimization, meaning
that inefficient perturbation schemes will be generated and that a
substantial amount of objective function evaluations will be wasted.
This is, for now, a critical weakness of the methodology. As future
work we intend to reduce the occurrence of those schemes by adding
learning mechanisms to DSO.

The DSO method introduced here is a naïve implementation, using random
distributions to generate scalings, to select methods for recombination
and bound correction, to select lower-quality solutions when stagnation
occurs, to update the firmware, among other aspects. Questions that
arise regarding the configuration and performance of the DSO algorithm
require future work, but the promising results indicate that the current
version is competitive with state-of-the-art methods.

\section*{Acknowledgments}

This paper was supported by the Brazilian Government CNPq (Universal)
grant (486950/2013-1) and CAPES (Science without Borders) grant (12180-13-0)
to V.V.M., and Canada\textquoteright s NSERC Discovery grant RGPIN
283304-2012 to W.B.

\bibliographystyle{IEEEtran}

\end{document}